\documentclass[journal]{ieeeconf}
\usepackage{lineno}
\usepackage[cmex10]{amsmath}
\usepackage{amscd,amsfonts,amsmath,amssymb,mathrsfs,pifont,stmaryrd,tipa}
\usepackage{array,cases,dsfont,graphicx,texdraw}
\usepackage{wrapfig}
\usepackage{graphics} 
\usepackage{epsfig} 
\usepackage{stfloats}
\usepackage{subfig}
\usepackage{hyperref}
\usepackage{color}
\newtheorem{Example}{Example}
\newtheorem{Remark}{Remark}
\newtheorem{Corollary}{Corollary}
\newtheorem{Definition}{Definition}

\newenvironment{Proof}{\noindent{\em Proof:\/}}{\hfill $\Box$\par}
\newtheorem{Theorem}{Theorem}
\newtheorem{Lemma}{Lemma}

\newtheorem{Assumption}{Assumption}

\input{tcilatex}                                                          

\IEEEoverridecommandlockouts                              
\title{\LARGE \bf
A Unified Strategy for Solution Seeking in Graphical $N$-coalition Noncooperative Games}

\author{Maojiao Ye, \emph{Member, IEEE}, Guoqiang Hu, \emph{Member, IEEE}, Frank L. Lewis, \emph{Fellow, IEEE}, and Lihua Xie, \emph{Fellow, IEEE}
\thanks{M. Ye is with the School of Automation, Nanjing University of Science and Technology, Nanjing 210094, P.R. China and was with the School of Electrical and Electronic Engineering, Nanyang Technological University, 639798 (Email: mjye@njust.edu.cn); G. Hu, L. Xie are with the School of Electrical and Electronic Engineering, Nanyang
Technological University, 639798, Singapore  (Email: gqhu@ntu.edu.sg, elhxie@ntu.edu.sg); F. L. Lewis is with the UTA Research Institute, University of Texas at Arlington, Fort Worth, Texas, USA (Email: lewis@uta.edu).}
\thanks{This work was supported by Singapore Economic Development Board under EIRP grant S14-1172-NRF EIRP-IHL, the National Natural Science Foundation of China (NSFC), No. 61803202, the Natural Science Foundation of Jiangsu Province, No. BK20180455 and the Fundamental Research Funds for the Central Universities, No. 30918011332.}
}

\begin{document}

\maketitle
\thispagestyle{empty}
\pagestyle{empty}

\begin{abstract}
This paper aims to reduce the communication and computation costs of the Nash equilibrium seeking strategy for the $N$-coalition noncooperative games proposed in \cite{YEAT}.  The objective is achieved in two manners: 1. An interference graph is introduced to describe the interactions among the agents in each coalition. 2. The Nash equilibrium seeking strategy is designed with the interference graphs considered. The convergence property of the proposed Nash equilibrium seeking strategy is analytically investigated. It is shown that the agents' actions generated by the proposed method converge to a neighborhood of the Nash equilibrium of the graphical $N$-coalition noncooperative games under certain conditions. Several special cases where there is only one coalition and/or there are coalitions with only one agent are considered. The results for the special cases demonstrate that the proposed seeking strategy achieves the solution seeking for noncooperative games, social cost minimization problems and single-agent optimization problems in a unified framework. Numerical examples are presented to support the theoretical results.
\end{abstract}

\vspace{-5mm}
\section{INTRODUCTION}

Recent years witnessed great efforts made by the researchers to study distributed optimization problems and noncooperative games (see e.g., \cite{YeTcyber17}-\cite{Ratliff16}). Distributed optimization captures the cooperative characteristics of networked systems and covers many practical applications including economic dispatch in the smart grids \cite{LiTII16} and resource allocation problems \cite{Xiao06}, etc. Noncooperative games, which catch the competitive nature of self-interested players, have been widely adopted to analyze electricity markets \cite{YeTcyber17}, transport systems \cite{Hollander06}, just to name a few. Motivated by the incentive to model the cooperative and competitive behaviors in economic markets and adversarial networked systems, an $N$-coalition noncooperative game was formulated in \cite{YEAT} to accommodate the cooperation and competition in a unified framework. The formulated game concerns with a set of agents that form $N$ interacting coalitions, each of which is composite of a subset of agents that cooperatively minimize the sum of their local objective functions. Nevertheless, the coalitions serve as self-interested players in noncooperative games.

The $N$-coalition noncooperative games can be utilized to accommodate many practical circumstances in which cooperation and competition coexist among the engaged participants. For example, in group competitions (e.g., football games), the participants in the same group cooperate to win the game while different groups are competitive opponents. From another aspect, agents in multi-agent environment (e.g., transportation networks, cloud computing) may benefit from forming coalitions to perform tasks especially when single agent could not complete the task independently or efficiently \cite{ShehoryAI98}. The agents in the same coalition work collaboratively to serve the tasks while distinct coalitions compete for the tasks to gain payoffs. Moreover, the formulated $N$-coalition noncooperative games have great potential to address networked systems in which the agents are subjected to multiple tasks. For instance, it is promising to achieve multi-party rendezvous where the agents within each coalition need to rendezvous to a certain location while preserving connectivity with the agents in other coalitions. Motivated by the practical relevance, this paper further investigates the $N$-coalition noncooperative games to reduce the communication and computation costs of the Nash equilibrium seeking strategy proposed in  \cite{YEAT}.

To relieve the communication and computation burden, a new Nash equilibrium seeking strategy is proposed for the $N$-coalition noncooperative games. The newly proposed seeking strategy follows the idea of implementing the gradient play by utilizing consensus protocols to estimate necessary information (see, e.g., \cite{YEAT}\cite{YeTcyber17}\cite{YeAC}). Consensus protocols were utilized as supportive tools to solve social cost minimization problems and noncooperative games in several recent works including \cite{YEAT}\cite{YeTcyber17}\cite{YeAC}-\cite{YECCC17} and some other references therein. In \cite{Nedich16}, the DIGing algorithm was proposed through imitating the distributed gradient descent via a gradient tracking recursion. A leader-following consensus protocol was leveraged to disseminate local information in \cite{YeAC} and \cite{YECCC17} for noncooperative games and distributed optimization problems, respectively. A dynamic average consensus protocol was included in an extremum seeking scheme in \cite{DoughertyTAC17}, \cite{YeTcyber17} and \cite{YEAT} for collaborative optimization, Nash equilibrium seeking for aggregative games and the $N$-coalition noncooperative games, respectively. However, in \cite{YEAT}, it assumed that the agents' objective functions are functions of all the agents' actions, which might possibly result in redundant communication and computation costs if the agents' objective functions depend on only a subset of the agents' actions.

Game theoretic models for wireless communication networks and multi-agent systems involve typical examples in which the players' payoff functions are more closely related with the actions' of the players who are their physical neighbors (see e.g., \cite{Chen12}-\cite{Marden09} and the references therein). Considering the spatial reuse and the heterogeneous  resource competition capabilities among the secondary users in the wireless communication networks, the asymmetric interference among the users was addressed in \cite{Chen12}, which generalized the symmetric interference in the atomic congestion games in \cite{Tekin12}. In the sensor deployment game designed in \cite{Marden09}, the local utilities for each node depended on the actions of their in-neighbors only. To describe the interactions among the agents, an interference graph was introduced \cite{SalehisadaghianiCDC14} \cite{Chen12}-\cite{YEICCA17}. The authors in \cite{Chen12} adopted a directed graph for the spatial spectrum access games and similar idea was utilized in \cite{Tekin12} to describe the interactions among the players in the atomic congestion games. In this regard,  the $N$-coalition game considered in \cite{YEAT} can be treated as a game in which the interference graphs for the coalitions are complete graphs.  To further remove the redundant estimation variables in the strategy design, we introduce an interference graph to each coalition and design an interference graph based Nash equilibrium seeking strategy, by which the communication and computation costs are further reduced if the interference graphs are not complete graphs.

In summary, with some preliminary results in \cite{YEICCA17}, the main contributions of the paper are listed as follows.
\begin{enumerate}
  \item A new Nash equilibrium seeking strategy is proposed for the $N$-coalition noncooperative games by introducing an interference graph to each coalition to describe the interactions among the agents in each coalition. Compared with the seeking strategy in \cite{YEAT}, the new strategy requires less communication and computation costs, especially for $N$-coalition noncooperative games on sparse interference graphs. Moreover, different from the seeking strategy in \cite{YEAT} which only addresses the case in which $N\geq 2, m_i\geq 2,\forall i\in\{1,2,\cdots,m_i\}$, the proposed seeking strategy can also accommodate the cases in which $N=1$ and/or $m_i=1.$ This indicates that the proposed method solves the single-agent minimization, social cost minimization and noncooperative games in a unified framework.
  \item The convergence results are analytically studied by Lyapunov stability analysis. It is shown that under certain conditions, the agents' actions converge to a neighborhood of the Nash equilibrium of the $N$-coalition noncooperative games for $N\geq 2$, and to a neighborhood of the solution of the corresponding minimization problem for $N=1.$
\end{enumerate}


\section{Problem Formulation}\label{p1_res}

 In this paper, we revisit the $N$-coalition noncooperative game formulated in \cite{YEAT}. In this game, there are $N$ interacting coalitions ($N$ is an integer and $N\geq 1$) that are self-interested to minimize their own cost functions. The cost function of coalition $i$, denoted as $f_i(\mathbf{x}_i,\mathbf{x}_{-i}),$ is defined as
\begin{equation}
f_i(\mathbf{x}_i,\mathbf{x}_{-i})=\sum\nolimits_{j=1}^{m_i} f_{ij}(\mathbf{x}_i,\mathbf{x}_{-i}), \forall i\in \{1,2,\cdots,N\},
\end{equation}
in which $m_i\geq 1$ is an integer that denotes the number of agents in coalition $i$ and $f_{ij}(\mathbf{x}_i,\mathbf{x}_{-i})$ is a function available to agent $j$ in coalition $i$ only. Furthermore,
$\mathbf{x}_i=[x_{i1},x_{i2},\cdots,x_{im_i}]^T,$
$\mathbf{x}_{-i}=[\mathbf{x}_1^T,\mathbf{x}_2^T,\cdots,\mathbf{x}_{i-1}^T,\mathbf{x}_{i+1}^T,\cdots,\mathbf{x}_N^T]^T,$
$\mathbf{x}=[\mathbf{x}_1^T,\mathbf{x}_2^T,\cdots,\mathbf{x}_N^T]^T,$
where $x_{ij}\in \mathbb{R}$ is the action of agent $j$ in coalition $i$ that is governed by agent $j$ in coalition $i$. Hence, in the $N$-coalition noncooperative game, the agents within the same coalition collaborate to minimize the sum of their local functions, while constituting a coalition that acts as a player in a noncooperative game.
The agents in coalition $i$ are equipped with a communication network with graph topology $\mathcal{G}_C^i,  \forall i\in \{1,2,\cdots,N\}$. This paper intends to develop a new Nash equilibrium seeking strategy for the agents to reduce the communication and computation costs of the strategy proposed in \cite{YEAT} under the condition that the Nash equilibrium of the $N$-coalition game $\mathbf{x}^*$ exists and is finite. Note that without any further clarification, if $N=1,$ the objective is solution seeking for the corresponding minimization problem given that the solution exists and is finite. Moreover, the Nash equilibrium considered in this paper is pure-strategy Nash equilibrium.

The Nash equilibrium of the $N$-coalition noncooperative game (for $N> 1$) is defined as follows.
\begin{Definition}
Nash equilibrium is an action profile on which no player can reduce its own cost by unilaterally changing its own action, i.e., $\mathbf{x}^*=(\mathbf{x}_i^*,\mathbf{x}_{-i}^*)$ is a Nash equilibrium of the $N$-coalition noncooperative games if for $\mathbf{x}_i\in \mathbb{R}^{m_i},$
\begin{equation}
f_i(\mathbf{x}_i^*,\mathbf{x}_{-i}^*)\leq f_i(\mathbf{x}_{i},\mathbf{x}_{-i}^*),\forall i\in \{1,2,\cdots,N\}.
\end{equation}
\end{Definition}
\begin{Remark}
Note that by regarding each coalition as a player in a noncooperative game, the conditions that ensure the existence and uniqueness of the Nash equilibrium follow the existing explorations in \cite{RosenEco65}\cite{Ratliff16}\cite{BasarPhi}-\cite{Facchinei07}.
\end{Remark}

The rest of the paper is based on the following assumptions.
\begin{Assumption}\label{smooth_assum}
The agents' objective functions $f_{ij}(\mathbf{x}_i,\mathbf{x}_{-i}),$ for all $i\in\{1,2,\cdots,N\},j\in\{1,2,\cdots,m_i\}$ are $\mathcal{C}^2$ functions.
\end{Assumption}
\begin{Assumption}\label{p_1_Assum_1}
For all distinct $\mathbf{x},\mathbf{y}\in \mathbb{R}^{\sum_{i=1}^Nm_i}$
\begin{equation}\label{monot}
(\mathbf{x}-\mathbf{y})^T\left(\mathcal{P}(\mathbf{x})-\mathcal{P}(\mathbf{y})\right)> 0,
\end{equation}
where $\mathcal{P}(\mathbf{x})=[(\frac{\partial f_1(\mathbf{x})}{\partial \mathbf{x}_1})^T,(\frac{\partial f_2(\mathbf{x})}{\partial \mathbf{x}_2})^T,\cdots,(\frac{\partial f_N(\mathbf{x})}{\partial \mathbf{x}_N})^T]^T.$
\end{Assumption}
\begin{Remark}
Assumption \ref{p_1_Assum_1}, adapted from \cite{YeAC} \cite{RosenEco65}-\cite{SalehisadaghianiCDC14},  assumes that the pseudo-gradient vector $\mathcal{P}(\mathbf{x})$ is strictly monotone. Note that this assumption weakens the strong monotonicity assumption (i.e., $(\mathbf{x}-\mathbf{y})^T(\mathcal{P}(\mathbf{x})-\mathcal{P}(\mathbf{y}))\geq l||\mathbf{x}-\mathbf{y}||^2,$ for some positive constant $l$) in \cite{YEAT}. Moreover, when $N=1,$ the strict (strong) monotonicity is reduced to strict (strong) convexity of the corresponding function.
Detailed explanations on monotonicity are available in \cite{Facchinei} and readers are referred to Chapter 2 of \cite{Facchinei} for more elaborations.
\end{Remark}

\section{$N$-coalition noncooperative games on interference graphs}\label{game_inter_graph}

In this section, a new Nash equilibrium seeking strategy will be proposed. Convergence results for the case in which $N\geq 2, m_i\geq 2,\forall i\in\{1,2,\cdots,N\}$ will be firstly investigated followed by several special cases in which $N=1$ and/or $m_i=1$ for some $i\in \{1,2,\cdots,N\}$.

Before we proceed to the Nash equilibrium seeking strategy design, we firstly introduce an interference graph (see e.g., \cite{SalehisadaghianiCDC14} \cite{Tekin12}) to each coalition to describe the interactions among the agents in the same coalition. Suppose that agent $j$ is a neighbor of agent $k$ in the interference graph of coalition $i$ (denoted as $\mathcal{G}_I^i$) if and only if $f_{ij}$ explicitly depends on $x_{ik}$ or alternatively $f_{ik}$ explicitly depends on $x_{ij}$ \cite{YEICCA17}.  Then,

\begin{equation}\label{equ_1}
\frac{\partial f_i(\mathbf{x})}{\partial x_{ik}}=\sum_{j=1}^{m_i}\frac{\partial f_{ij}(\mathbf{x})}{\partial x_{ik}}=\sum_{j\in (\mathcal{N}_{Ik}^i\cup \{k\})}\frac{\partial  f_{ij}(\mathbf{x})}{\partial x_{ik}},
\end{equation}
where $\mathcal{N}_{Ik}^i$ denotes the neighboring set of agent $k$ in coalition $i$ in the interference graph $\mathcal{G}_I^i.$
\begin{figure}
\begin{center}
\scalebox{0.5}{\includegraphics[74,341][587,560]{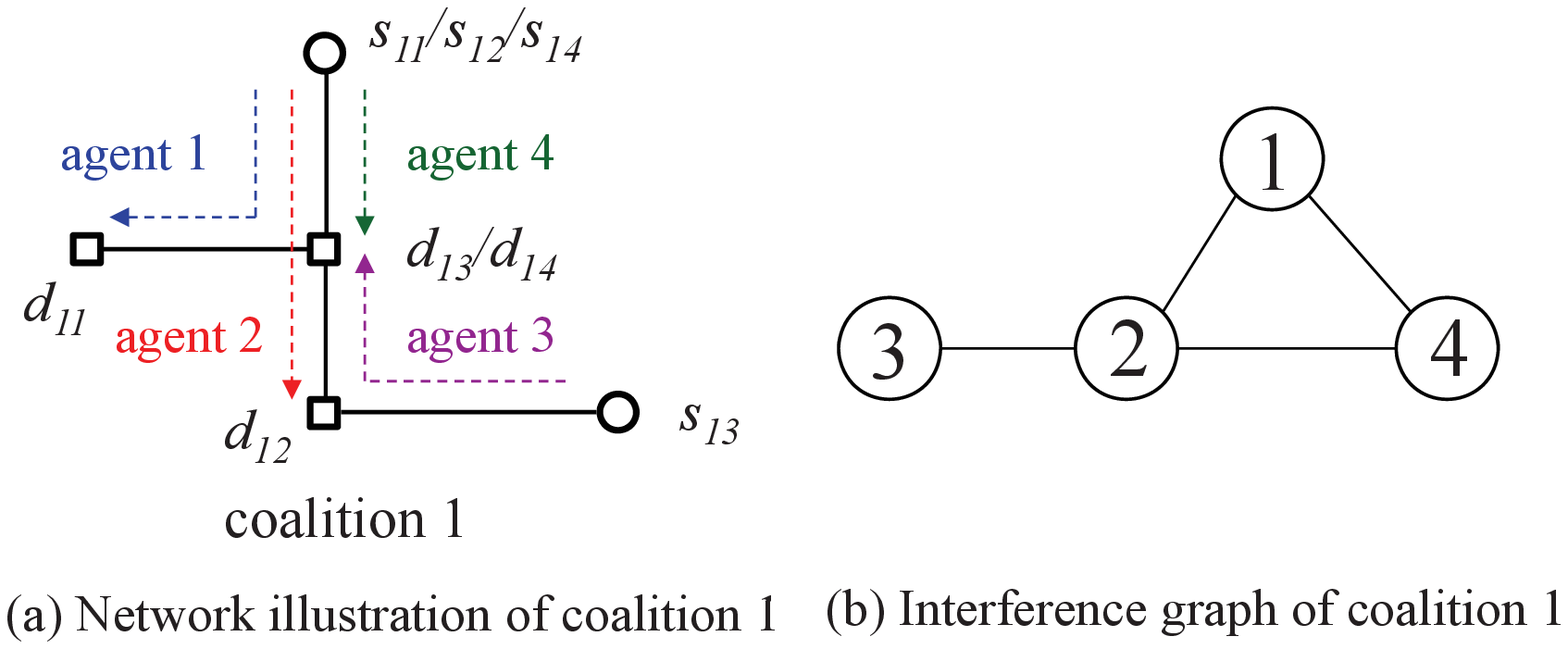}}
\caption{Network illustration and the interference graph for the agents(users) in coalition $1$ of the game in Example \ref{Examp_coal1}.}\label{inter_c1}
\end{center}
\end{figure}

\begin{Example}\label{Examp_coal1}
Consider an $N$-coalition network game in Ad Hoc wireless communication networks. In the game, each agent $j$ in coalition $i$ intends to send its data flow $x_{ij}$ from the source node $s_{ij}$ to the destination node $d_{ij}$. Suppose that coalition $1$ of the $N$-coalition noncooperative game contains agents $1$-$4$, whose associated network is depicted in Fig. \ref{inter_c1}(a) (note that the agents in other coalitions are not depicted in the figure). Then, by the congestion control model in \cite{Alpan02}, the cost functions for the $4$ agents are $f_{11}(x_{11},x_{12},x_{14},\mathbf{x}_{-1}),f_{12}(x_{11},x_{12},x_{13}$,
$x_{14},\mathbf{x}_{-1}),f_{13}(x_{12},x_{13},$$\mathbf{x}_{-1})$ and $f_{14}(x_{11},x_{12},x_{14},\mathbf{x}_{-1})$, respectively. In this example, the interference graph for coalition $1$ can be depicted as in Fig. \ref{inter_c1}(b).
In addition, $\frac{\partial f_1(\mathbf{x})}{\partial x_{11}}=\frac{\partial f_{11}(\mathbf{x})}{\partial x_{11}}+\frac{\partial f_{12}(\mathbf{x})}{\partial x_{11}}+\frac{\partial f_{14}(\mathbf{x})}{\partial x_{11}},$ $\frac{\partial f_{1}(\mathbf{x})}{\partial x_{12}}=\frac{\partial f_{11}(\mathbf{x})}{\partial x_{12}}+\frac{\partial f_{12}(\mathbf{x})}{\partial x_{12}}+\frac{\partial f_{13}(\mathbf{x})}{\partial x_{12}}+\frac{\partial f_{14}(\mathbf{x})}{\partial x_{12}},$ $\frac{\partial f_{1}(\mathbf{x})}{\partial x_{13}}=\frac{\partial f_{12}(\mathbf{x})}{\partial x_{13}}+\frac{\partial f_{13}(\mathbf{x})}{\partial x_{13}},$ $\frac{\partial f_{1}(\mathbf{x})}{\partial x_{14}}=\frac{\partial f_{11}(\mathbf{x})}{\partial x_{14}}+\frac{\partial f_{12}(\mathbf{x})}{\partial x_{14}}+\frac{\partial f_{14}(\mathbf{x})}{\partial x_{14}},$ which are in line with \eqref{equ_1}.
\end{Example}

\begin{Remark}
Equation \eqref{equ_1} is inspired by practical engineering systems in which the agents' costs are more closely correlated with their physical neighbors (see e.g., \cite{Chen12}\cite{Tekin12}\cite{Alpan02}).
By introducing the interference graph to capture the player-specific cost functions in such games, we intend to propose a Nash equilibrium seeking strategy for the $N$-coalition noncooperative games with the interactions among the agents taken into account.
\end{Remark}


\subsection{Nash equilibrium seeking strategy design}

To search for the Nash equilibrium of the $N$-coalition noncooperative games, agent $j$ in coalition $i$ updates its action according to
\begin{equation}\label{strate_1_1}
\dot{x}_{ij}(t)=-d_{ij}g_{ijj}(t),
\end{equation}
for $i\in \{1,2,\cdots,N\},j\in \{1,2,\cdots,m_i\}$ where $d_{ij}=\delta \bar{d}_{ij}$, $\delta$ is a small positive parameter to be determined and $\bar{d}_{ij}$ is a fixed positive parameter. Moreover, for each
$i\in\{1,2,\cdots,N\},j\in\{1,2,\cdots,m_i\},k\in (\mathcal{N}_{Ij}^i \cup \{j\}),$ $g_{ijk}$ is generated by
\begin{equation}\label{strate_1_2}
\begin{aligned}
\dot{w}_{ijk}(t)&=-\sum\nolimits_{l\in (\mathcal{N}_{Ik}^i\cup \{k\})}a_i^{jl}(g_{ijk}(t)-g_{ilk}(t))\\
g_{ijk}(t)&=w_{ijk}(t)+\frac{\partial f_{ij}(\mathbf{x}(t))}{\partial x_{ik}},
\end{aligned}
\end{equation}
where $a_i^{jl}$ is the element on the $j$th row and $l$th column of the adjacency matrix of $\mathcal{G}_C^i.$ In addition, $w_{ijk}$ for $i\in\{1,2,\cdots,N\},j\in\{1,2,\cdots,m_i\},k\in (\mathcal{N}_{Ij}^i \cup \{j\})$ are auxiliary variables initialized at $w_{ijk}(0)=0$ (without loss of generality, the initial time instant is supposed to be $0$). Note that in the subsequent analysis, the time-variable $t$ might be omitted for notational convenience.


\begin{Remark}
The design of the Nash equilibrium seeking strategy for the $N$-coalition noncooperative game in \eqref{strate_1_1}-\eqref{strate_1_2} follows the ideas in \cite{YEAT}\cite{YeAC} to adapt the consensus protocols in \cite{Freeman06} (i.e., \eqref{strate_1_2} in the proposed method) to estimate the required information based on which the gradient play is indirectly implemented (i.e., \eqref{strate_1_1} in the proposed method). To be more specific, $g_{ijj}$ provides an estimate of $\frac{1}{(N_{Ij}^i+1)}\frac{\partial f_{i}(\mathbf{x})}{\partial x_{ij}}$, where $N_{Ij}^i$ denotes the number of agents in $\mathcal{N}_{Ij}^i$, by utilizing the physical-interaction-based consensus protocol in \eqref{strate_1_2}. Moreover, \eqref{strate_1_1} implements the gradient play by utilizing the agents' local estimates on the corresponding gradient values. However, it's worth noting that in \eqref{strate_1_1}-\eqref{strate_1_2}, agent $j$ in coalition $i$ only generates $2(N_{Ij}^i+1)$ auxiliary variables (i.e., $g_{ijk}$,$w_{ijk},k\in (\mathcal{N}_{Ij}^i\cup \{j\})$). Moreover, agent $j$ in coalition $i$ only communicates $g_{ijk}$ for $k\in (\mathcal{N}_{Ij}^i\cup \{j\})$ with the agents who are in $\mathcal{N}_{Ik}^i\cup\{k\}$ and are its neighbors in $\mathcal{G}_C^i$. In contrast, in \cite{YEAT}, agent $j$ in coalition $i$ needs to generate $2m_i$ auxiliary variables and communicate all the auxiliary variables with its neighbors in the communication graph $\mathcal{G}_C^i$. Hence, when the interference graphs are sparse graphs, $N_{Ij}^i+1$ is much smaller compared to $m_i$ and the method in \eqref{strate_1_1}-\eqref{strate_1_2} gives much lower computation and communication costs compared with the method in \cite{YEAT}. Furthermore, even if the interference graphs are complete graphs $(N_{Ij}^i+1=m_i),i\in\{1,2,\cdots,N\},j\in\{1,2,\cdots,m_i\}$, agent $j$ in coalition $i$ only communicates with its neighbors in $\mathcal{G}_C^i$ on $g_{ijk}$ but not on $w_{ijk},$ where $k\in\{1,2,\cdots,m_i\}$, which is also required in \cite{YEAT}.
\end{Remark}

\subsection{Convergence analysis}
 In the following, we firstly provide the analysis for the case in which $N\geq 2, m_i\geq 2,\forall i\in \{1,2,\cdots,N\}.$ Then, the special cases where $m_i=1$ and/or $N=1$ will be discussed.

\emph{Case $\uppercase\expandafter{\romannumeral1}$ ($N\geq 2,m_i\geq 2,\forall i\in \{1,2,\cdots,N\}$):}
For each $i\in\{1,2,\cdots,N\},k\in\{1,2,\cdots,m_i\},$ remove from $\mathcal{G}_C^i$ the nodes that are not the neighbors of agent $k$ in coalition $i$ in the interference graph $\mathcal{G}_I^i$ except for agent $k$. Leave the remaining nodes and the edges therein \cite{YEICCA17}. Denote the new graph as $\mathcal{G}_{C_k}^i$, named as interference-to-$k$ communication graph for coalition $i$,  and its corresponding Laplacian matrix as $L_{C_k}^i.$  Moreover, define $W_{ik}/G_{ik}/\bar{F}_{ik}(\mathbf{x})$ as the concatenated vectors of $w_{ijk}/g_{ijk}/\frac{\partial f_{ij}(\mathbf{x})}{\partial x_{ik}}$ over $j\in (\mathcal{N}_{Ik}^i\cup \{k\}),$ respectively. Then, \eqref{strate_1_2} can be written as
\begin{equation}\label{con_ve}
\dot{W}_{ik}=-L_{C_k}^iG_{ik},G_{ik}=W_{ik}+\bar{F}_{ik}(\mathbf{x}).
\end{equation}
Let $\mathbf{1}_{M}(\mathbf{0}_M)$ be an $M$-dimensional column vector composite of $1(0)$ and define  $N_{C_k^i}$ as the number of agents in $\mathcal{G}_{C_k}^i$. Motivated by \cite{YeTcyber17}, let $R_{ik}\in \mathbb{R}^{N_{C_k^i}}\times \mathbb{R}^{N_{C_k^i}-1}$ be a matrix such that $U_{ik}=\left[\frac{\mathbf{1}_{N_{C_k^i}}}{\sqrt{N_{C_k^i}}} \ \ R_{ik}\right]\in \mathbb{R}^{N_{C_k^i}}\times \mathbb{R}^{N_{C_k^i}}$ is an orthogonal matrix. Note that $\left[\frac{\mathbf{1}_{N_{C_k^i}}}{\sqrt{N_{C_k^i}}} \ \ R_{ik}\right]$  stands for a matrix whose first column is  $\frac{\mathbf{1}_{N_{C_k^i}}}{\sqrt{N_{C_k^i}}}$ and the remaining columns are equal to those in $R_{ik}.$
In addition, define
$[\underline{G}_{ik} \ \ \bar{G}_{ik}^T]^T=U_{ik}^T G_{ik},$ $[\underline{W}_{ik}\ \ \bar{W}_{ik}^T]^T=U_{ik}^T W_{ik},$ where $\underline{G}_{ik},\underline{W}_{ik}\in \mathbb{R}$ and  $\bar{G}_{ik},\bar{W}_{ik}\in \mathbb{R}^{N_{C_k^i}-1}.$

Then, given that $\mathcal{G}_{C_k}^i$ is undirected and connected,

\begin{equation}
\dot{\underline{W}}_{ik}=0,\dot{\bar{W}}_{ik}=-R_{ik}^TL_{C_k}^i R_{ik} \bar{G}_{ik},
\end{equation}
by which $\underline{W}_{ik}(t)=\underline{W}_{ik}(0)=0,$ and $R_{ik}^TL_{C_k}^iR_{ik}$ is symmetric positive definite.

Since $G_{ik}=W_{ik}+\bar{F}_{ik}(\mathbf{x}),$
\begin{equation}
\underline{G}_{ik}=\underline{W}_{ik}+\frac{\mathbf{1}_{N_{C_k^i}}^T}{\sqrt{N_{C_k^i}}}\bar{F}_{ik}(\mathbf{x}),\bar{G}_{ik}=\bar{W}_{ik}+R_{ik}^T\bar{F}_{ik}(\mathbf{x}),
\end{equation}
i.e., $\underline{G}_{ik}(t)=\underline{W}_{ik}(t)+\frac{1}{\sqrt{N_{C_k^i}}}\sum_{j\in (\mathcal{N}_{Ik}^i\cup \{k\})}\frac{\partial f_{ij}(\mathbf{x})}{\partial x_{ik}}=\frac{1}{\sqrt{N_{C_k^i}}}\sum_{j=1}^{m_i}\frac{\partial f_{ij}(\mathbf{x})}{\partial x_{ik}}$ by \eqref{equ_1}, given that $\mathcal{G}_{C_k}^i$ is undirected and connected. In addition,
\begin{equation}
\dot{\bar{G}}_{ik}=-R_{ik}^TL_{C_k}^iR_{ik}\bar{G}_{ik}+R_{ik}^T\dot{\bar{F}}_{ik}(\mathbf{x}).
\end{equation}

\begin{Lemma}\label{bound_1}
For $\mathbf{x}\in \mathbb{R}^{\sum_{i=1}^Nm_i},$ there exists a positive constant $\beta_{ijk}$ such that
\begin{equation}
\left|\left|g_{ijk}-\frac{1}{N_{C_k^i}}\sum_{l=1}^{m_i}\frac{ \partial f_{il}(\mathbf{x})}{\partial x_{ik}}\right|\right|\leq \beta_{ijk}||\bar{G}_{ik}||,
\end{equation}
for each $i\in \{1,2,\cdots,N\},k\in \{1,2,\cdots,m_i\},j\in (N_{Ik}^i\cup \{k\})$ given that  for each $i\in \{1,2,\cdots,N\},k\in \{1,2,\cdots,m_i\},$ $\mathcal{G}_{C_k}^i$ is undirected and connected.
\end{Lemma}
\begin{Proof}
See Section \ref{bound_1_prof} for the proof.
\end{Proof}

Let $\bar{G}=[\bar{G}_{11}^T,\bar{G}_{12}^T,\cdots,\bar{G}_{1m_1}^T,\bar{G}_{21}^T,\cdots,\bar{G}_{2m_2}^T,\cdots,$ $\bar{G}_{Nm_N}^T]^T.$ Define $\chi(t)=[\bar{G}^T(t),(\mathbf{x}(t)-\mathbf{x}^*)^T]^T.$ Then, the following result can be derived.

\begin{Theorem}\label{int_res1}
Suppose that for each $i\in\{1,2,\cdots,N\},k\in \{1,2,\cdots,m_i\},$ the interference graph $\mathcal{G}_I^i$ and the interference-to-$k$ communication graphs $\mathcal{G}_{C_k}^i$ are undirected and connected.  Then, for each pair of positive constants $(\Delta,v)$, there is a positive constant $\delta^*(\Delta,v)$ such that for each $\delta\in (0,\delta^*),$ $\chi(t)$ generated by \eqref{strate_1_1}-\eqref{strate_1_2} satisfies
\begin{equation}
||\chi(t)||\leq \phi(||\chi(0)||,\delta t)+v
\end{equation}
where $\phi(\cdot,\cdot)\in\mathcal{KL}$ and $||\chi(0)||<\Delta,$ given that Assumptions \ref{smooth_assum}-\ref{p_1_Assum_1} hold.
\end{Theorem}
\begin{Proof}
See Section \ref{int_res1_proof} for the proof.
\end{Proof}

%

\emph{Case $\uppercase\expandafter{\romannumeral2}$ ($N\geq 2,m_j=1$ for some $j\in \{1,2,\cdots,N\}$):} For notational convenience, we suppose that $m_i>1,\forall i\neq j$, which can be easily extended to the case in which there are multiple coalitions with only one agent.
In this case, the seeking strategy in \eqref{strate_1_1}-\eqref{strate_1_2} for the agent in coalition $j$ is reduced to
\begin{equation}\label{reduced_game}
\dot{x}_{j1}=-d_{j1}g_{j11},\dot{w}_{j11}=0,g_{j11}=w_{j11}+\frac{\partial f_{j1}(\mathbf{x})}{\partial x_{j1}},
\end{equation}
where $w_{j11}(0)=0, d_{j1}=\delta \bar{d}_{j1}.$

Hence, \eqref{reduced_game} is equivalent to
\begin{equation}\label{int_red}
\dot{x}_{j1}=-d_{j1}\frac{\partial f_{j1}(\mathbf{x})}{\partial x_{j1}}.
\end{equation}

Define $\bar{G}_{-j}$ as the concatenated vector of $\bar{G}_{ik}$ for $i\in \{1,2,\cdots,N\},i\neq j,k\in \{1,2,\cdots,m_i\}$ and $\chi_{-j}(t)=[\bar{G}_{-j}^T(t),(\mathbf{x}(t)-\mathbf{x}^*)^T]^T.$ Then, the following result can be derived.
\begin{Corollary}\label{p2_res_1_col1}
Suppose that for each $i\in\{1,2,\cdots,N\},i\neq j, k\in \{1,2,\cdots,m_i\},$ the interference graph $\mathcal{G}_I^i$ and the interference-to-$k$ communication graphs $\mathcal{G}_{C_k}^i$ are undirected and connected. Then, for each pair of positive constants $(\Delta,v)$, there exists a positive constant $\delta^*(\Delta,v)$ such that for each $\delta\in (0,\delta^*)$, $\chi_{-j}(t)$ generated by \eqref{strate_1_1}-\eqref{strate_1_2} satisfies
\begin{equation}
||\chi_{-j}(t)||\leq \phi(||\chi_{-j}(0)||,\delta t)+v,
\end{equation}
where $\phi(\cdot,\cdot)\in \mathcal{KL}$ and $||\chi_{-j}(0)||<\Delta,$ given that Assumptions \ref{smooth_assum}-\ref{p_1_Assum_1} hold.
\end{Corollary}
\begin{Proof}
See Section \ref{p1_res_1_col1_proof} for the proof.
\end{Proof}

In addition, if $m_i=1,\forall i\in \{1,2,\cdots,N\}$ and $N\geq 2$, the considered $N$-coalition noncooperative game is reduced to a \emph{\textbf{noncooperative game}} in which there exist $N$ players and the cost function of player $i$ is $f_{i1}(\mathbf{x}).$ Moreover, the proposed seeking strategy in \eqref{strate_1_1}-\eqref{strate_1_2} is reduced to
 \begin{equation}\label{reduce_grad}
\dot{x}_{i1}=-d_{i1}\frac{\partial f_{i1}(\mathbf{x})}{\partial x_{i1}}, \forall i\in\{1,2,\cdots,N\} ,
\end{equation}
by which it can be derived that the Nash equilibrium of the noncooperative game is globally asymptotically stable given that Assumptions \ref{smooth_assum}-\ref{p_1_Assum_1} hold.


\emph{Case $\uppercase\expandafter{\romannumeral3}$ ($N=1,m_1\geq 2$):} In this case,  there is only one coalition and there are multiple agents therein. The  formulated $N$-coalition game is reduced to a \emph{\textbf{social cost minimization problem}}, i.e.,
\begin{equation}\label{soc_min}
\text{min}_{\mathbf{x}_1} \ \ \ f_1(\mathbf{x}_1)=\sum\nolimits_{j=1}^{m_1}f_{1j}(\mathbf{x}_1).
\end{equation}
Moreover, the seeking strategy in \eqref{strate_1_1}-\eqref{strate_1_2} is reduced to
\begin{equation}\label{strate_r_1}
\begin{aligned}
&\dot{x}_{1j}=-d_{1j}g_{1jj}\\
&\dot{w}_{1jk}=-\sum\nolimits_{l\in (\mathcal{N}_{Ik}^1\cup \{k\})}a_1^{jl}(g_{1jk}-g_{1lk})\\
&g_{1jk}=w_{1jk}+\frac{\partial f_{1j}(\mathbf{x}_1)}{\partial x_{1k}},
\end{aligned}
\end{equation}
in which $d_{1j}=\delta \bar{d}_{1j}, w_{1jk}(0)=0,$ $\forall j\in \{1,2,\cdots,m_1\},k\in (\mathcal{N}_{Ij}^1\cup\{j\}).$ Define $\chi_1(t)=[\bar{G}_1^T(t),(\mathbf{x}_1(t)-\mathbf{x}_1^*)^T]^T,$  where $\mathbf{x}_1^*$ is the solution to the problem in \eqref{soc_min} and $\bar{G}_1=[\bar{G}_{11}^T,\bar{G}_{12}^T,\cdots,\bar{G}_{1m_1}^T]^T$. Then, following the analysis of Theorem \ref{int_res1}, the subsequent corollary can be obtained.

\begin{Corollary}\label{p1_res_1_col2}
Suppose that Assumptions \ref{smooth_assum}-\ref{p_1_Assum_1} hold, $\mathcal{G}_I^1$ and $\mathcal{G}_{C_k}^1,\forall k\in\{1,2,\cdots,m_1\}$ are undirected and connected. Then, for each pair of positive constants $(\Delta,v)$, there exists a positive constant $\delta^*(\Delta,v)$ such that for each $\delta\in (0,\delta^*)$, $\chi_1(t)$ generated by \eqref{strate_r_1} satisfies
\begin{equation}
||\chi_1(t)||\leq \phi(||\chi_1(0)||,\delta t)+v,
\end{equation}
where $\phi(\cdot,\cdot)\in\mathcal{KL}$ and  $||\chi_1(0)||<\Delta$, given that Assumptions \ref{smooth_assum}-\ref{p_1_Assum_1} hold.
\end{Corollary}

\begin{Remark}
In this case, Assumption \ref{p_1_Assum_1} is reduced to the condition that $f_1(\mathbf{x}_1)$ is strictly convex and Corollary \ref{p1_res_1_col2} indicates that social cost minimization can be achieved without convexity conditions on the local objective functions.
\end{Remark}

\emph{Case $\uppercase\expandafter{\romannumeral4}$ ($N=1,m_1=1$):}
In this case, the formulated $N$-coalition game is reduced to a \emph{\textbf{single-agent minimization problem}}, i.e.,
\begin{equation}\label{min_pr1}
\text{min}_{x_{11}}  \ \ f_{11}(x_{11}),
\end{equation}
and the proposed seeking strategy is reduced to

\begin{equation}\label{desc_ref}
\dot{x}_{11}=-d_{11}\frac{\partial f_{11}(x_{11})}{\partial x_{11}},
\end{equation}
by which it can be derived that $x_{11}^*,$ which stands for the solution to the problem in \eqref{min_pr1}, is globally asymptotically stable under \eqref{desc_ref} given that Assumptions \ref{smooth_assum}-\ref{p_1_Assum_1} hold.


\begin{Remark}
The above results indicate that the considered $N$-coalition noncooperative game covers the social cost minimization, noncooperative games and single-agent optimization as special cases. Moreover, the proposed seeking strategy solves the aforementioned problems in a unified framework.
\end{Remark}

\begin{Remark}
In the above results, we suppose that for each $i\in \{1,2,\cdots,N\},k\in\{1,2,\cdots,m_i\},$ $\mathcal{G}_{C_k}^i$ is undirected and connected. When the interference graphs are complete graphs, this condition is reduced to $\mathcal{G}_C^i$ for each $i\in\{1,2,\cdots,N\}$ being undirected and connected. With the interference graphs, the connectivity of $\mathcal{G}_{C_k}^i$ for $i\in\{1,2,\cdots,N\},k\in \{1,2,\cdots,m_i\}$ is utilized to ensure that the agents can disseminate necessary information to the other related agents. The following condition, adopted from \cite{SalehisadaghianiCDC14}, is utilized to ensure that $\mathcal{G}_{C_k}^i$ for $i\in \{1,2,\cdots,N\},k\in \{1,2,\cdots,m_i\}$ are connected.
\end{Remark}

\begin{Assumption}\label{inter_gra}
\cite{SalehisadaghianiCDC14} The interference graphs $\mathcal{G}_I^i$ and the communication graphs $\mathcal{G}_C^i$ for all $i\in \{1,2,\cdots,N\}$ are undirected and connected. Moreover,
$\mathcal{G}_m^i \subseteq \mathcal{G}_C^i\subseteq \mathcal{G}_I^i,\forall i\in \{1,2,\cdots,N\},$ where $\mathcal{G}_m^i$ denotes the maximal triangle-free spanning subgraph of $\mathcal{G}_I^i$.
\end{Assumption}

Then, following the analysis in \cite{SalehisadaghianiCDC14}, the following result can be derived.
\begin{Lemma}\label{gra_1}
Suppose that Assumption \ref{inter_gra} is satisfied. Then, for each $i\in\{1,2,\cdots,N\},k\in \{1,2,\cdots,m_i\},$ the graph $\mathcal{G}_{C_k}^i$ is undirected and connected.
\end{Lemma}

\begin{Remark}
In the above results, we suppose that the interference graphs for all the coalitions with more than one agent are connected. However, if there is an interference graph $\mathcal{G}_I^j$ for some $j\in \{1,2,\cdots,N\}$ that is not connected, then, the coalition can be separated into several coalitions (possibly into several independent agents), each with a connected interference graph. For more detailed elaborations on the interference graphs, the reader are referred to \cite{YEICCA17}.
\end{Remark}

\begin{Remark}
Note that with strong monotonicity of the pseudo-gradient vector, it can be shown that similar to the results in \cite{YEAT}, the proposed method enables the agents' actions to converge exponentially to the Nash equilibrium under certain conditions.
\end{Remark}

\section{Numerical Examples}\label{P_1_exa}
In this section, $N$-coalition noncooperative games with $N=3,m_1=1,m_2=3,m_3=6$ are considered. Firstly, a congestion control game in Ad Hoc wireless communication networks is considered in Example \ref{examp1}. Then, a game that violates Assumption \ref{p_1_Assum_1} is simulated in Example \ref{examp2}.

\begin{figure}
\begin{center}
\scalebox{0.35}{\includegraphics[69,257][594,539]{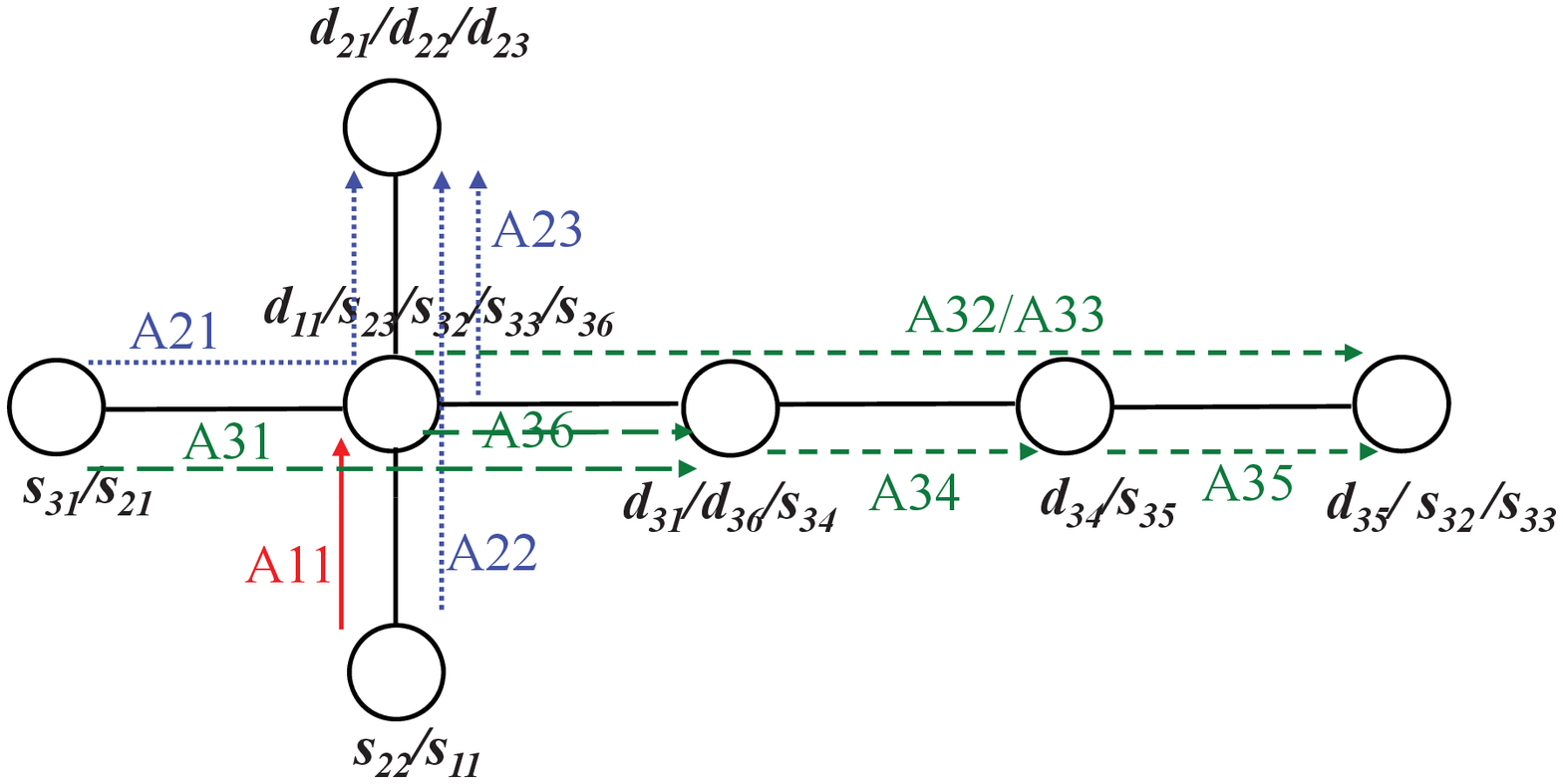}}
\caption{Network illustration for Example \ref{examp1}.}\label{source2}
\end{center}
\end{figure}

\begin{Example}\label{examp1}
Consider an Ad Hoc wireless communication network that consists of several nodes as depicted in Fig. \ref{source2}. Suppose that in the network, agent $j$ in coalition $i$ intends to transmit nonnegative data flow $x_{ij}$ from $s_{ij}$ to $d_{ij}$. Then, the cost function of agent $j$ in coalition $i$ can be captured by  \cite{Alpan02}
\begin{equation}
f_{ij}(\mathbf{x})=\sum_{l\in \hat{P}_{ij}}\frac{\kappa}{C_l-\sum_{l}x_{ij}}-u_{ij}\log(x_{ij}+1),
\end{equation}
where $C_l$ is the capacity of link $l$, $\kappa,$ $u_{ij}$ are agent-specified parameters, $\hat{P}_{ij}$ denotes the path of agent $j$ in coalition $i$, and $\sum_{l}x_{ij}$ denotes the total flow on link $l,l\in \hat{P}_{ij}$. See \cite{Alpan02} for the details on the model formulation of the congestion control game.

In the simulation, we suppose that all the link capacities are 20 and $\kappa=10,u_{ij}=10$. By direct calculation, it can be derived that the game admits a unique Nash equilibrium $\mathbf{x}^*=[12.63,5.58,3.68,6.12,5.16,2.03,2.03,10.16,10.16,5.63]^T.$
\end{Example}

\begin{figure}
\begin{center}
\scalebox{0.31}{\includegraphics[100,323][470,535]{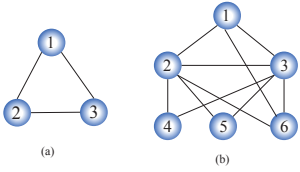}}
\caption{The interference graphs for coalitions $2$ and $3$. (a) is the interference graph for the agents in coalition $2$ and (b) is the interference graph for the agents in coalition $3$.}\label{interf_graph}
\end{center}
\end{figure}

\begin{figure}
\begin{center}
\scalebox{0.35}{\includegraphics[82,328][455,525]{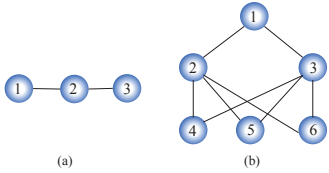}}
\caption{The communication graphs for the agents in coalitions $2$ and $3$. (a) is the communication graph for the agents in coalition $2$ and (b) is the communication graph for the agents in coalition $3$.}\label{comm_graph_2}
\end{center}
\end{figure}

\begin{figure}
\begin{center}
\scalebox{0.43}{\includegraphics[84,285][482,525]{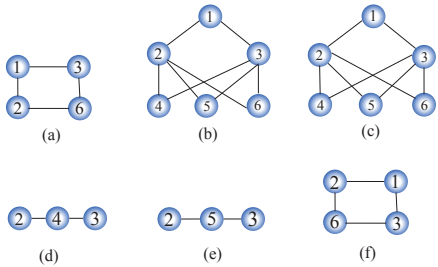}}
\caption{The graphs $\mathcal{G}^3_{C_k}$ for $k\in\{1,2,\cdots,6\}$ are given in (a)-(f), respectively.}\label{reduce_graph}
\end{center}
\end{figure}

\begin{figure}
\begin{center}
\scalebox{0.4}{\hspace{-20mm}\includegraphics[0,0][397,273]{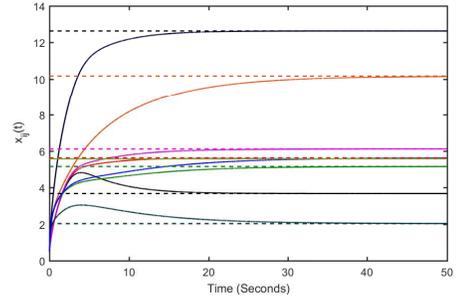}}
\caption{The agents' actions generated by the proposed method in \eqref{strate_1_1}-\eqref{strate_1_2} in Example \ref{examp1}.}\label{action_inf1}
\end{center}
\end{figure}

The interference graph and the communication graph for the game are given in Figs. \ref{interf_graph}-\ref{comm_graph_2}, respectively. Correspondingly, $\mathcal{G}^3_{C_k}$ for $k\in\{1,2,\cdots,6\}$ are given in Fig. \ref{reduce_graph}. It can be seen that in line with Lemma \ref{gra_1}, $\mathcal{G}^3_{C_k}$ for $k\in\{1,2,\cdots,6\}$ are undirected and connected.
Initialized at $x_{ij}(0)=0.5$ (all the other variables are initialized at zero), the players' actions generated by the proposed method in  \eqref{strate_1_1}-\eqref{strate_1_2}
are given in Fig. \ref{action_inf1}.
The simulation result shows that the agents' actions generated by the proposed seeking strategy in \eqref{strate_1_1}-\eqref{strate_1_2} converge to the actual Nash equilibrium.

Moreover, different from \cite{YEAT}, by designing the strategy based on the interference graphs in coalition $3,$ agent $1$ does not need to generate $g_{314},w_{314},g_{315},w_{315}$; agent $4$ does not need to generate $g_{341},w_{341},g_{345},w_{345},g_{346},w_{346};$ agent $5$ does not need to generate $g_{351},w_{351},g_{354},w_{354},g_{356},w_{356}$ and agent $6$ does not need to generate $g_{364},w_{364},g_{365},w_{365}$ in the seeking strategy, thus requiring less computation. In addition, from the communication graph for the agents in coalition $3$, it can be seen that if the method in \cite{YEAT} is adopted, agent $1$ in coalition 3 needs to communicate with agents $2$ and $3$ in coalition $3$ on $g_{31k},w_{31k}$; agent $2$ in coalition $3$ needs to communicate with agents $1,4,5,6$ in coalition $3$ on $g_{32k},w_{32k}$; agent $3$ in coalition $3$ needs to communicate with agents $1,4,5,6$ in coalition $3$ on $g_{33k},w_{33k}$; agents $4,5,6$ in  coalition $3$ need to communicate with agents $2$ and $3$ on $g_{34k},w_{34k}$, $g_{35k},w_{35k}$ and $g_{36k},w_{36k}$, $k\in\{1,2,\cdots,6\},$ respectively. In contrast, by utilizing the proposed method, agent $1$ in coalition $3$ communicates with agents $2,3$ in coalition $3$ on $g_{31k},k\in\{1,2,3,6\}$; agents $2,3$ in coalition $3$ communicates with agents $1,4,5,6$ on $g_{32k},g_{33k},k\in \{1,2,\cdots,6\},$ respectively; agents $4,5,6$ in coalition $3$ communicate with agents $2,3$ on $g_{34k},k\in\{2,3,4\},g_{35k},k\in\{2,3,5\},g_{36k},k\in\{1,2,3,6\},$ respectively. Through direct comparison, it can be seen that by utilizing the proposed method, the communication and computation costs are significantly reduced compared with the method in \cite{YEAT}.

\begin{Example}\label{examp2}
\textbf{(A game without strictly monotone pseudo-gradients)} In this example,
we consider a game in which the agents' objective functions are defined as $f_{11}(\mathbf{x})=(x_{11}-\frac{1}{2}x_{31})^2,$ $f_{21}=x_{21}^3+2x_{21}^2+x_{22}^2-x_{21}x_{32}-x_{21}x_{23},$
$f_{22}=-2x_{21}x_{22}+x_{21}x_{23},$
$f_{23}=-x_{21}^3+x_{23}^2-x_{21}x_{11}+x_{21}x_{22},$ $f_{31}=-e^{x_{31}}+x_{31}^2-x_{31}x^2_{21},$ $f_{32}=-x_{32}\sum_{j=1}^6 x_{3j},$ $f_{33}=x_{32}^2+x_{33}^2+x_{32}\sum_{j=1}^6 x_{3j}-\sum_{i=1}^2\sum_{j=1}^{m_i}x_{ij},$
$f_{34}=x_{34}^2,$ $f_{35}=x_{35}^2,$ and $f_{36}=e^{x_{31}}+x_{36}^2,$ respectively. In this case, letting $\mathcal{P}(\mathbf{x})=\mathbf{0}_{10}$ gives two points on which $\mathbf{x}=\mathbf{0}_{10}$ or $\mathbf{x}=[49,14,7, 0,98,0,0,0,0,0]^T,$ by which it can be easily derived that Assumption
\ref{p_1_Assum_1} is violated.
\end{Example}

\begin{figure}
\begin{center}
\scalebox{0.48}{\includegraphics[127,257][474,522]{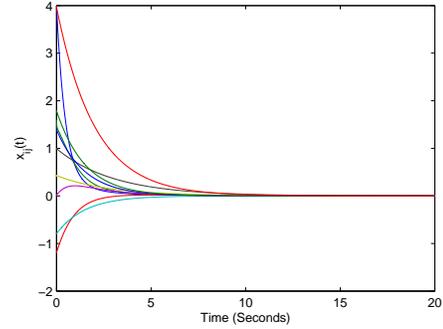}}
\caption{The agents' actions generated by the proposed method in \eqref{strate_1_1}-\eqref{strate_1_2} in Example \ref{examp2}.}\label{action_inf2}
\end{center}
\end{figure}

With the agents' actions initialized at $\mathbf{x}(0)=[4,1.6,-1.2,-0.8,0,0.4,1,1.4,1.8,4]^T$ (all the other variables are initialized at zero), the agents' actions generated by the proposed method are depicted in Fig. \ref{action_inf2}. Note that the communication graphs and the interference graphs adopted in the simulation are the same as those in Example \ref{examp1}.
From the simulation result, we can see that the agents' actions generated by the proposed method converge to $\mathbf{0}_{10}$ under the given initial conditions, which is a Nash equilibrium point of the game, though Assumption \ref{p_1_Assum_1} is not satisfied. Actually, similar to \cite{YEAT}, local convergence to each Nash equilibrium that satisfies certain conditions (see Assumptions 3-4 in \cite{YEAT}) can be derived for the proposed seeking strategy. The analytical study is not provided due to space limitation.

\section{Conclusion}\label{P_1_CONC}

This paper investigates the $N$-coalition noncooperative games formulated in \cite{YEAT}. To reduce the computation and communication costs of the Nash equilibrium seeking strategy proposed in \cite{YEAT}, a new seeking strategy is proposed by introducing an interference graph to each coalition to describe the interactions among the agents. A sufficient condition on the interference graph and the communication graph is provided to ensure the convergence of the proposed seeking strategy. It is shown that under certain conditions, the agents' actions generated by the proposed methods converge to a neighborhood of the Nash equilibrium of the $N$-coalition noncooperative games. In addition, several special cases, where there is only one coalition and/or there are coalitions with only one agent, are studied. The associated convergence results demonstrate that the proposed method solves the noncooperative games, the social cost minimization problems and the single-agent optimization problems in a unified framework.

\section{Appendix}
\subsection{Preliminaries}\label{preli_coal}
A graph is defined as $\mathcal{G}=(\mathcal{V},\mathcal{E})$ where $\mathcal{V}=\{1,2,\cdots,M\}$ is the set of nodes in the
network and $\mathcal{E}$ denotes the set of edges. The elements of $\mathcal{E}$ are denoted as $(i,j)$, which represents for an edge from $i$ to $j.$ Associate with each edge $(j,i)\in \mathcal{E}$ a weight $a_{ij}>0.$ Note that we suppose that the graph is simple, i.e., $a_{ii}=0.$ The graph is undirected if $a_{ij}=a_{ji},\forall i,j\in \mathcal{V}.$ The adjacency matrix $\mathcal{A}$ of graph $\mathcal{G}$ is a matrix with its element on the $i$th row and $j$th column defined as
$a_{ij}>0$ if $(j,i)\in \mathcal{E},$ else, $a_{ij}=0.$ An undirected graph is connected if there exists a path between any pair of distinct vertices.  The neighboring set of agent $i$ is defined as $\mathcal{N}_i=\{j \in \mathcal{V}|(j,i)\in \mathcal{E}\}.$
The Laplacian matrix $L$ is defined as $L=\mathcal{D}-\mathcal{A}, $ where $\mathcal{D}$ is a diagonal matrix whose $i$th diagonal element is equal to $\sum_{j=1}^{M}a_{ij}$.

A subgraph $H_{\mathcal{G}}$ of a graph $\mathcal{G}$ is a graph whose sets of vertices and edges are subsets of the vertex set and edge set of $\mathcal{G}$, respectively. Alternatively, $\mathcal{G}$ is the supergraph of $H_{\mathcal{G}}.$ A subgraph $H_{\mathcal{G}}$ is a spanning subgraph of $\mathcal{G}$ if it contains all the vertices of $\mathcal{G}.$ In addition, a triangle-free subgraph $H_{\mathcal{G}}$ of a graph $\mathcal{G}$ is a subgraph in which no three vertices form a triangle of edges. Graph $H_{\mathcal{G}}$ is a maximal triangle-free subgraph of $G$ if adding a missing edge to $H_{\mathcal{G}}$ forms a triangle \cite{SalehisadaghianiCDC14}\cite{Godsil2013}.
\subsection{Proof of Lemma \ref{bound_1}}\label{bound_1_prof}

For notational convenience, given that $g_{ijk}$ is the $l$th element in $G_{ik}$, we use $[R_{ik}\bar{G}_{ik}]_l$ to denote the $l$th element in $R_{ik}\bar{G}_{ik}.$
Since $
G_{ik}=\left[\frac{\mathbf{1}_{N_{C_{k}^i}}}{\sqrt{N_{C_{k}^i}}} \ \ R_{ik}\right][\underline{G}_{ik},(\bar{G}_{ik})^T]^T=\mathbf{1}_{N_{C_{k}^i}} \left(\frac{1}{N_{C_{k}^i}}\sum_{j=1}^{m_i}\frac{\partial f_{ij}(\mathbf{x})}{\partial x_{ik}}\right)+R_{ik} \bar{G}_{ik},$ it can be derived that $\left|\left|g_{ijk}-\frac{1}{N_{C_{k}^i}}\sum_{l=1}^{m_i}\frac{ \partial f_{il}(\mathbf{x})}{\partial x_{ik}}\right|\right|
=\left|\left|\frac{1}{N_{C_{k}^i}}\sum_{l=1}^{m_i}\frac{ \partial f_{il}(\mathbf{x})}{\partial x_{ik}}+[R_{ik}\bar{G}_{ik}]_{l}-\frac{1}{N_{C_{k}^i}}\sum_{l=1}^{m_i}\frac{ \partial f_{il}(\mathbf{x})}{\partial x_{ik}}\right|\right|
= \left|\left|[R_{ik}\bar{G}_{ik}]_{l}\right|\right|\leq \beta_{ijk}||\bar{G}_{ik}||,
$
where $\beta_{ijk}$ is a positive constant such that $||[R_{ik}\bar{G}_{ik}]_{l}||\leq \beta_{ijk} ||\bar{G}_{ik}||.$

\subsection{Proof of Theorem \ref{int_res1}}\label{int_res1_proof}

Inspired by \cite{YEAT}\cite{YeAC} and \cite{Khailil}, define the Lyapunov candidate function as $V=\sum_{i=1}^N \sum_{k=1}^{m_i} \bar{G}_{ik}^TP_{ik}\bar{G}_{ik}+\frac{1}{2}(\mathbf{x}-\mathbf{x}^*)^T\left(\text{diag}\{\frac{\bar{d}_{ij}}{N_{C_j^i}}\}\right)^{-1}(\mathbf{x}-\mathbf{x}^*),
$
where $P_{ik}$ is a symmetric positive definite matrix such that $P_{ik}R_{ik}^TL_{C_k}^iR_{ik}+R_{ik}^TL_{C_k}^iR_{ik}P_{ik}=Q_{ik},$ $Q_{ik}$ is a symmetric positive definite matrix. Moreover, the notation $\text{diag}\{h_{ij}\}$ for $i\in \{1,2,\cdots,N\},j\in\{1,2,\cdots,m_i\}$ is a diagonal matrix whose diagonal elements are $h_{11},h_{12},\cdots,h_{1m_1},h_{21},\cdots, h_{Nm_N},$ successively. 
%

Define $\lambda_{min}(Q_{ik})$ as the minimum eigenvalue of $Q_{ik},$ and $\underline{\lambda}=\text{min}\{\lambda_{min}(Q_{ik}),i\in \{1,2,\cdots,N\},k\in\{1,2,\cdots,m_i\}\}.$ Then, the time derivative of the Lyapunov candidate function satisfies
\begin{equation}
\begin{aligned}
\dot{V}&\leq -\underline{\lambda}||\bar{G}||^2-\delta(\mathbf{x}-\mathbf{x}^*)^T\mathcal{P}(\mathbf{x})\\
&+\delta (\mathbf{x}-\mathbf{x}^*)^T\left(\mathcal{P}(\mathbf{x})-\text{col}\{N_{C_j^i}g_{ijj}\}\right)\\
&-2\delta \sum_{i=1}^N\sum_{k=1}^{m_i} \bar{G}_{ik}^TP_{ik}R_{ik}^T\left(\frac{\partial \bar{F}_{ik}(\mathbf{x})}{\partial \mathbf{x}}\right)^T\text{col}\{\bar{d}_{ij}g_{ijj}\},
\end{aligned}
\end{equation}
where $\text{col}\{h_{ij}\}$ for $i\in\{1,2,\cdots,N\},j\in\{1,2,\cdots,m_i\}$ is defined as $\text{col}\{h_{ij}\}=[h_{11},h_{12},\cdots,h_{1m_1},h_{21},\cdots,h_{Nm_N}]^T.$ Note that the above inequality is derived by utilizing the fact that
$(\mathbf{x}-\mathbf{x}^*)^T \text{col}\{N_{C_j^i}g_{ijj}\}=(\mathbf{x}-\mathbf{x}^*)^T (\mathcal{P}(\mathbf{x})-\mathcal{P}(\mathbf{x})+\text{col}\{N_{C_j^i}g_{ijj}\})$
and $\dot{\bar{F}}_{ik}=\left(\frac{\partial \bar{F}_{ik}(\mathbf{x})}{\partial \mathbf{x}}\right)^T\dot{\mathbf{x}}=-\delta\left(\frac{\partial \bar{F}_{ik}(\mathbf{x})}{\partial \mathbf{x}}\right)^T \text{col}\{\bar{d}_{ij}g_{ijj}\}$.

By Lemma \ref{bound_1}, it can be derived that for $\mathbf{x}\in \mathbb{R}^{\sum_{i=1}^N m_i}$, there is a positive constant $\beta$ such that $\left|\left| \mathcal{P}(\mathbf{x})-\text{col}\{N_{C_j^i}g_{ijj}\}\right|\right|\leq \beta ||\bar{G}||$. In addition,  as $f_{ij}(\mathbf{x})$ for $i\in\{1,2,\cdots,N\},j\in \{1,2,\cdots,m_i\}$ are $\mathcal{C}^2$ functions, $||\frac{\partial \bar{F}_{ik}(\mathbf{x})}{\partial \mathbf{x}}||$ is bounded for $\mathbf{x}$ that belongs to a compact set. From another aspect, by Assumption \ref{p_1_Assum_1}, $\mathcal{P}(\mathbf{x}^*)=\mathbf{0},$ as $x_{ij}\in \mathbb{R}$ in this paper. Hence, there are positive constants $\beta_1,\beta_2$ such that
$||\text{col}\{N_{C_j^i}g_{ijj}\}||=||\text{col}\{N_{C_j^i}g_{ijj}\}-\mathcal{P}(\mathbf{x})+\mathcal{P}(\mathbf{x})-\mathcal{P}(\mathbf{x}^*)||\leq \beta_1||\bar{G}||+\beta_2||\mathbf{x}-\mathbf{x}^*||$ for $\mathbf{x}$ that belongs to a compact set.

Since $(\mathbf{x}-\mathbf{x}^*)^T\mathcal{P}(\mathbf{x})=(\mathbf{x}-\mathbf{x}^*)^T(\mathcal{P}(\mathbf{x})-\mathcal{P}(\mathbf{x}^*))>0$ for all $\mathbf{x}\neq \mathbf{x}^*$, there exists a function $\varphi(||\mathbf{x}-\mathbf{x}^*||)\in \mathcal{K}$ such that $(\mathbf{x}-\mathbf{x}^*)^T\mathcal{P}(\mathbf{x})\geq \varphi(||\mathbf{x}-\mathbf{x}^*||)$ by Assumption \ref{p_1_Assum_1}. Therefore, for any $\chi=[\bar{G}^T,(\mathbf{x}-\mathbf{x}^*)^T]^T$ that belongs to a compact set $\Omega$ that contains the origin, there are positive constants $\gamma_1,\gamma_2$ such that
\begin{equation}
\begin{aligned}
\dot{V}\leq &-\underline{\lambda}||\bar{G}||^2 -\delta\varphi(||\mathbf{x}-\mathbf{x}^*||)+\delta \beta || \mathbf{x}-\mathbf{x}^*||||\bar{G}||\\
&+\delta ||\bar{G}||(\gamma_1||\bar{G}||+\gamma_2||\mathbf{x}-\mathbf{x}^*||)\\
\leq & -\left(\underline{\lambda}-\delta \gamma_1-\frac{\delta(\beta+\gamma_2)}{2\gamma_3}\right)||\bar{G}||^2-\delta \varphi(||\mathbf{x}-\mathbf{x}^*||)\\
&+\delta (\beta+\gamma_2)\gamma_3||\mathbf{x}-\mathbf{x}^*||^2/2,
\end{aligned}
\end{equation}
where $\gamma_3$ is a positive constant that can be arbitrarily chosen.

Suppose that $\underline\lambda-\delta \gamma_1-\frac{\delta(\beta+\gamma_2)}{2\gamma_3}>\delta.$ Then,
\begin{equation}
\begin{aligned}
\dot{V}\leq &-\frac{\delta}{2}(||\bar{G}||^2+\varphi(||\mathbf{x}-\mathbf{x}^*||))-\frac{\delta}{2}(||\bar{G}||^2+\varphi(||\mathbf{x}-\mathbf{x}^*||))\\
&+\delta(\beta+\gamma_2)\gamma_3||\mathbf{x}-\mathbf{x}^*||^2/2\leq  -\delta \varphi_1(||\chi||)/2,
\end{aligned}
\end{equation}
for all $\varphi_1(||\chi||)\geq (\beta+\gamma_2)\gamma_3 \varrho,$ where $\varrho$ is a positive constant such that $||\mathbf{x}-\mathbf{x}^*||^2\leq \varrho$ for $\chi\in \Omega$ and $\varphi_1\in\mathcal{K}$  satisfies $||\bar{G}||^2+\varphi(||\mathbf{x}-\mathbf{x}^*||)\geq \varphi_1(||\chi||).$

For $\chi\in \Omega,$ choose $\gamma_3$ such  that $(\beta+\gamma_2)\gamma_3\varrho$ is sufficiently small. Then, fixed $\gamma_3$, let $\delta^*=\frac{2\underline{\lambda}\gamma_3}{2\gamma_3+\beta+\gamma_2+2\gamma_1\gamma_3}.$

Then, for each $\delta\in (0,\delta^*)$
\begin{equation}
\dot{V}\leq -\frac{\delta}{2}\varphi_1(||\chi||),\forall ||\chi||\geq \varphi_1^{-1}((\beta+\gamma_2)\gamma_3\varrho).
\end{equation}
By further noticing that there are positive constants $c_1$ are $c_2$ such that $c_1||\chi||^2\leq V \leq c_2||\chi||^2,$ the conclusion can be derived by following the proof of Theorem 4.18 in  \cite{Khailil}.

\subsection{Proof of Corollary \ref{p2_res_1_col1}}\label{p1_res_1_col1_proof}

The result can be proven by defining the Lyapunov candidate function as $V=\sum_{i=1,i\neq j}^N\sum_{k=1}^{m_i}\bar{G}_{ik}^TP_{ik}\bar{G}_{ik}+\frac{1}{2}(\mathbf{x}-\mathbf{x}^*)^T\left(\text{diag}\left\{\frac{\bar{d}_{il}}{N_{C_{l}^i}}\right\}\right)^{-1}(\mathbf{x}-\mathbf{x}^*).
$
The rest of the proof is similar to that of Theorem \ref{int_res1}, and is omitted here.

\end{document}